%


%
%

\def\today{\ifcase\month\or January\or February\or
March\or April\or May\or June\or July\or August\or
September\or October\or November\or December\fi
\space\number\day, \number\year}




\def\dspace{\lineskip=2pt\baselineskip=18pt
\lineskiplimit=0pt}
\def\sspace{\lineskip=2pt\baselineskip=12pt
\lineskiplimit=0pt}

\font \bbrm=cmbx10 at 12pt

\font \ninerm= cmr10 at 9pt

\def\smalltype{\ninerm}

\def\bigtype{\bbrm}

\hsize=13.5cm
\magnification=1200
\def\ce{\centerline}

\def\hb{\hfill\break}

\def\title #1{\null\bigskip\ce{\bigtype #1}
\bigskip}

\def\alp{\alpha}		
\def\bet{\beta}		
\def\gam{\gamma}		
\def\del{\delta}

\def\tet{\theta}		\def\Tet{\Theta}

\def\kap{\kappa}

\def\ome{\omega}		


\def\calF{{\cal F}}

\def\calK{{\cal K}}

\def\calP{{\cal P}}



    
\font\tenboldgreek=cmmib10
 \font\sevenboldgreek=cmmib10 at 7pt
\font\fiveboldgreek=cmmib10 at 7pt
\newfam\bgfam
\textfont\bgfam=\tenboldgreek
\scriptfont\bgfam=\sevenboldgreek
\scriptscriptfont\bgfam=\fiveboldgreek

\mathchardef\ggarrow="7010

\font\tengerman=eufm10 \font\sevengerman=eufm7
\font\fivegerman=eufm5
\font\tendouble=msym10 \font\sevendouble=msym7
\font\fivedouble=msym5

\textfont4=\tengerman \scriptfont4=\sevengerman
\scriptscriptfont4=\fivegerman
\newfam\dbfam
\textfont\dbfam=\tendouble \scriptfont\dbfam=
\sevendouble
\scriptscriptfont\dbfam=\fivedouble

\mathchardef\ng="702D
\mathchardef\dbA="7041
\mathchardef\sm="7072
\mathchardef\nvdash="7030
\mathchardef\nldash="7031
\mathchardef\lne="7008
\mathchardef\sneq="7024
\mathchardef\spneq="7025
\mathchardef\sne="7028
\mathchardef\spne="7029
\mathchardef\ltms="706E
\mathchardef\tmsl="706F

\mathchardef\dbA="7041


\mathchardef\dbA="7041 
\mathchardef\dbB="7042 
\mathchardef\dbC="7043 \def\CC{{\fam=\dbfam\dbC}}
\mathchardef\dbD="7044 
\mathchardef\dbE="7045 
\mathchardef\dbF="7046 
\mathchardef\dbG="7047 
\mathchardef\dbH="7048 
\mathchardef\dbI="7049 
\mathchardef\dbJ="704A 
\mathchardef\dbK="704B 
\mathchardef\dbL="704C 
\mathchardef\dbM="704D 
\mathchardef\dbN="704E 
\mathchardef\dbO="704F 
\mathchardef\dbP="7050 
\mathchardef\dbQ="7051 
\mathchardef\dbR="7052 
\mathchardef\dbS="7053 
\mathchardef\dbT="7054 
\mathchardef\dbU="7055 
\mathchardef\dbV="7056 
\mathchardef\dbW="7057 
\mathchardef\dbX="7058 
\mathchardef\dbY="7059 
\mathchardef\dbZ="705A 

\def\nek{,\ldots,}
\def\sdp{\times \hskip -0.3em {\raise 0.3ex
\hbox{$\scriptscriptstyle |$}}} 


\def\min{\mathop{\rm min}}










\def\ddownarrow{\big\downarrow \hskip-0.70em\raise
2pt\hbox {$\big\downarrow$}}
\def\longright #1#2 {\smash{\mathop{\hbox to
#1pt {\rightarrowfill}}\limits_{#2}}}
\def\sqr#1#2{{\vcenter{\hrule height.#2pt\hbox{\vrule
width.#2pt height#1pt \kern#1pt \vrule width.#2pt}
\hrule height.#2pt}}}

\def\buildrul#1\under#2{\mathrel{\mathop{\null#2}
\limits_{#1}}}

\def\boxit#1{\vbox{\hrule\hbox{\vrule\kern3pt
\vbox{\kern3pt#1 \kern3pt}\kern3pt\vrule}\hrule}}

\def\subheading#1{\medskip\goodbreak\noindent{\bf
#1.}\quad}

\def\sect#1{\goodbreak\bigskip\centerline{\bf#1}
\medskip}
\def\pr{\smallskip\noindent{\bf Proof:\quad}}
\def\onumber #1{\ooalign{\hfil\raise.07ex\hbox{
\hfill$\scriptstyle \,#1$\hfil}
\cr\cr{$\bigcirc$}}}
\def\onumber c{\ooalign{\hfil\raise.07ex\hbox
{\hfill$\scriptstyle \,c$\hfil}
\cr\cr{$\bigcirc$}}}
\def\alpcirc {\ooalign{\hfil\raise.07ex
\hbox{\hfill$\scriptstyle\alp\;$\hfill}\cr\cr
{$\bigcirc$}}}
\def\astcirc {\ooalign{\hfil\raise.07ex
\hbox{\hfill$\textstyle\ast\;$\hfill}\cr\cr
{$\bigcirc$}}}

\def\longmapright #1#2 {\smash{\mathop{\hbox to
#1pt {\rightarrowfill}}\limits^{#2}}}
\def\longmapleft #1 #2 {\smash{\mathop{\hbox to
#1 pt {\leftarrowfill}}\limits^{#2}}}

\def\references#1{\goodbreak\bigskip\par\centerline
{\bf References}\medskip\parindent=#1pt}
\def\ref#1{\par\smallskip\hang\indent\llap{\hbox
to \parindent{#1\hfil\enspace}}\ignorespaces}

\def\back{{\raise 2.5pt\hbox{$\,\scriptscriptstyle
\backslash\,$}}}

\def\part{\partial}
\def\lwr #1{\lower 5pt\hbox{$#1$}\hskip -3pt}
\def\rse #1{\hskip -3pt\raise 5pt\hbox{$#1$}}
\def\lwrs #1{\lower 4pt\hbox{$\scriptstyle #1$}
\hskip -2pt}
\def\rses #1{\hskip -2pt\raise 3pt\hbox
{$\scriptstyle #1$}}

\def\<#1{\left\langle{#1}\right\rangle}

\def\subinbn{{\subset\hskip-8pt\raise 0.95pt
\hbox{$\scriptscriptstyle\subset$}}}

\def\llvdash{\mathop{\|\hskip-2pt
\raise 3pt\hbox{\vrule height 0.25pt width 1.5cm}}}

\def\lvdash{\mathop{|\hskip-2pt \raise 3pt\hbox
{\vrule height 0.25pt width 1.5cm}}}

\def\fakebold#1{\leavevmode\setbox0=\hbox{#1}%
  \kern-.025em\copy0 \kern-\wd0
  \kern .025em\copy0 \kern-\wd0
  \kern-.025em\raise.0333em\box0 }

\font\msxmten=msxm10
\font\msxmseven=msxm7
\font\msxmfive=msxm5
\newfam\myfam
\textfont\myfam=\msxmten
\scriptfont\myfam=\msxmseven
\scriptscriptfont\myfam=\msxmfive
\mathchardef\rhookupone="7016
\mathchardef\ldh="700D
\mathchardef\leg="7053
\mathchardef\ANG="705E
\mathchardef\lcu="7070
\mathchardef\rcu="7071
\mathchardef\leseq="7035
\mathchardef\qeeg="703D
\mathchardef\qeel="7036
\mathchardef\blackbox="7004
\mathchardef\bbx="7003
\mathchardef\simsucc="7025

\def\rhookup{{\fam=\myfam \rhookupone}}

\def\bigsquare{{\fam=\myfam\bbx}}

\font\tencaps=cmcsc10
\def\smallcaps{\tencaps}

\def\author#1{\bigskip\ce{\smallcaps #1}\medskip}

\def\upddots{\mathinner{\mkern
1mu\raise 1pt \hbox{.}\mkern 2mu \mkern
2mu \raise 4pt\hbox{.}\mkern 1mu \raise 7pt\vbox
{\kern 7 pt\hbox{.}}} }

\def\varchi{\ooalign{{\raise
1.385pt\hbox{$\chi$}}\crcr\hbox{--}\crcr}}

\def\trianarrow{{\raise 2pt\hbox to 0.50cm
{\hrulefill}\triangleright}}
\def\Chi{{\raise 3pt\hbox{$\chi$}}}

\font\b=cmr10 scaled \magstep4

\def\bigzerou{\smash{\lower1.7ex\hbox{\b 0}}}
\def\bigast{\smash{\lower1.7ex\hbox{\b *}}}

\def\leaderfill{leaders\hbox to 5em{\hss\hss}\hfill}
\newcount\notenumber
       
       \def\note#1{\advance\notenumber by
1\footnote{$^{\the\notenumber}$} {\sspace\smalltype #1}}

\null
\overfullrule=0pt
\def\llvdash{\mathop{\|\hskip-1.2pt \raise 3pt\hbox{\vrule
height 0.25pt width .4cm}}}
\def\Cub{{\rm Cub}}
{\nopagenumbers
\sect{SOME RESULTS ON NONSTATIONARY IDEAL
II}
\bigskip

\ce{Moti Gitik}
\medskip
\ce{School of Mathematical Sciences}
\ce{Sackler Faculty of Exact Sciences}
\ce{Tel Aviv University}
\ce{Tel Aviv 69978, Israel}
\dspace
\vskip 1truecm
\sect{Abstract}

{\medskip\narrower 

\noindent
We answer some question of [Gi].  The
upper bound of [Gi] on the strength of
$NS_{\mu^+}$ precipitous for a regular $\mu$ is proved
to be exact.  It is shown that
saturatedness of $NS_\kap^{\aleph_0}$ over
inaccessible $\kap$ requires at least $o(\kap) =
\kap^{++}$.  The upper bounds on the
strength of $NS_\kap$ precipitous for
inaccessible $\kap$ are reduced quite close to
the lower bounds.\medskip}
\vfill\eject
}
\dspace
\sect{0.~~Introduction}

The paper is a continuation of [Gi].  An
understanding of [Gi] is required.
However, there is one exception, Proposition
2.1.  It does not require any previous
knowledge and we think it is interesting on
its own.

The paper is organized as follows: In Section
1 we examine the strength of $NS_{\mu^+}$
precipitous.  The proof of the main theorem
there is a continuation of the proof of
2.5.1 from [Gi].  Section 2 deals with
saturation and answers question 3 of [Gi].
In Section 3 a new forcing construction of
$NS_\kappa$ precipitous over inaccessible
is sketched.  It combines ideas from [Gi,
Sec. 3] and [Gi1].  We assume  familiarity
with these papers.

We are grateful to the referee for his
remarks and suggestions.

\sect{1.~~On the strength of
precipitousness over a successor of
regular}

Our aim will be to improve the results of
[Gi] on precipitousness of $NS_{\mu^+}$ for
regular $\mu$ to the equiconsistency. 
Throughout the paper $\calK (\calF)$ is the
Mitchell Core Model with the maximal
sequence of measures $\calF$, under the
assumption $(\neg \exists \alp \ o^\calF(\alp) =
\alp^{++}$).  $o^\calF(\kap)$ denotes the
Mitchell order of $\kap$ or in other words
the length of the sequence $\calF$ over
$\kap$.  We refer to Mitchell [Mi1] for precise
definitions.

In order to state the result let us recall
a notion of $(\ome, \del)$-repeat point
introduced in [Gi].

\subheading{Definition}  Let $\alp, \del$
be ordinals with $\del < o^\calF(\kap)$.
Then $\alp$ is called a $(\ome,
\del)$-repeat point if (1)  $cf\alp =
\ome$, (2) for every $A \in \cap \{
\calF(\kap, \alp') | \alp \le \alp' < \alp
+ \del\}$ there are unboundedly many $\gam$'s
in $\alp$ such that $A \in \cap \{ \calF
(\kap, \gam') | \gam \le \gam' < \gam +
\del \}$.

We are going to prove the following:

\proclaim Theorem 1.1.  Suppose $NS_{\mu^+}$ is
precipitous for a regular $\mu>\aleph_1$ and
$GCH$.  Then there exists an $(\ome, \mu +
1)$--repeat point over $\mu^+$ in $K(\calF)$.

\subheading{Remark}

It is shown in [Gi] that starting with an
$(\ome, \mu + 1)$--repeat point it is
possible to obtain a model of $NS_{\mu^+}$
precipitous.  On the other hand
precipitousness of $NS_{\mu^+}^{\aleph_0}$
implies $(\ome,\mu)$--repeat point. 

In what follows we will actually continue
the proof of 2.5.1 of [Gi] and assuming that the 
$NS_{\mu^+}$ is precipitous (or even only
$NS_{\mu^+}^{\aleph_0}$ and
$NS_{\mu^+}^\mu$) we will obtain
$(\ome, \mu + 1)$--repeat point.

\pr 
Let $\kap = \mu^+$. 
We consider the ordinal $\alp^* < o^\calF(\kap)$
of the proof of 2.5.1 [Gi].  It was shown
there to be a $(\ome, \mu)$-repeat point,
under the assumption of nonexistence of
up-repeat point.  Intuitively, one can consider
 $\alp^*$ as the least relevant ordinal.
Basically, an ordinal $\alp$ is
called relevant if some condition in
$NS_\kap$ forces that the measure $\calF(\kap,
\alp)$ is used first in the generic
ultrapower to move $\kap$ and the
cofinality of $\kap$ changes to $\ome$.
Using a nonexistence of up-repeat point, a
set $A \in \calF(\kap, \alp^*)$ such that
$A \not\in \calF(\kap, \bet)$ for $\bet,
o^\calF(\kap) > \bet > \alp^*$, was picked.
This set $A$ was used in [Gi] and will be
used here to pin down $\alp^*$. Thus, for
$\tau \le \kap$ if there exists a largest
$\tau_1 < o^ \calF(\bet)$ such that $A
\cap \tau \in \calF(\tau, \tau_1)$ then we
denote it by $\tau^*$.  In this notation
$\kap^*$ is just $\alp^*$.  If $E = \{ \tau
< \kap |$ there exists $\tau^* \}$ then $E
\in \calF(\kap, \bet)$ for every $\bet$
with $\alp^* < \bet < o^\calF(\kap)$.  Also,
$A \cup E$ contains all points of
cofinality $\ome$ of a club, since by the
definition of $\alp^*$, $A \cup E \in
\cap \{ \calF(\kap, \alp) | \alp$ is a relevant
ordinal$\}$.

\proclaim Claim 1.  The set of $\alp< \kap$
satisfying (a) and (b) below is stationary
in $\kap$.
\item{(a)} $cf\ \alp = \mu$;
\item{(b)} for every $i < \mu$
$$\{\bet < \alp \mid cf\ \bet = \aleph_0\ {\rm
and}\ o^\calF(\bet) \ge \bet^* + i\}$$
is a stationary subset of $\alp$. 

\pr Otherwise, let $C$ be a club avoiding
all the $\alp$'s which satisfy (a) and (b).
Let $N$ be a good model in the sense of 2.5.1
of [Gi], with $C \in N$.  Consider $\langle
\tau_n^N | n < \ome\rangle$, $\langle d^N_n | n <
\ome\rangle$ and $\langle \bet^*_n | n < \ome \rangle$ of
2.5.1 [Gi].  Recall that $\langle \tau_n^N | n <
\ome \rangle$ is a sequence of indiscernibles for
$N$, each $\tau_n^N$ is a limit point of
$C$, $d_n^N$ is an $\ome$-club in $\cup (N
\cap \tau_n)$ consisting of indiscernibles
of cofinality $\ome$ in $C$, for $\nu \in
d_n^N$\enskip $\nu^*$ exists and $\bet_n^*$
represents it over $\kap$ (identically for
every $\nu, \nu' \in d_n^N)$. Also for
every $\tau < \tau'$ in $d_n^N$\enskip $\bet^N
(\tau) < \bet^N (\tau')$, where $\bet^N
(\tau)$ is the index of the measure over
$\kap$ to which $\tau$ corresponds. 

Fix $n < \ome$.  Then, $\tau_n \in C$.  As
in 2.10 or 2.14 of [Gi] we can assume that $cf
\tau_n = \mu$.  Since  (b) fails,
there are $i_n < \mu$ and $C_n$ a club of
$\tau_n$ disjoint with $\{ \nu < \tau_n |
cf \nu = \aleph_0$ and $o^{\calF} (\nu)
\ge \nu^* + i_n\}$.  Using elementarity of
$N$, it is easy to find such $C_n$ inside
$N$. Let $\del = \bigcup\limits_{n < \ome} i_n$.
Using 2.1.1 (or 2.15 for inaccessible
$\mu$) of [Gi] we will obtain $N^*
\supseteq N$ which agrees (mod initial
segment) with $N$ about indiscernibles but
has sets $d_n^{N^*}$ long enough to reach
$\del$, i.e. there will be a final segment
of $\tau$'s in $d_n^{N^*}$ with $\bet^{N^*}
(\tau) > \bet^*_n + \del$.  But then, for
such $\tau$, $o^\calF(\tau) \ge \tau^* +
\del$.  This is impossible, since $C_n$,
$d_n^{N^*}$ are both clubs of $\tau_n$ in
$N^*$ with bounded intersection.
Contradiction.
\hb
$\bigsquare$

Let $S$ denote the set of $\alp$'s
satisfying the conditions (a) and (b) of
Claim 1.  Now form a generic ultrapower
with $S$ in the generic ultrafilter.
Denote it by $M$ and let $\calF(\kap, \xi)$
be the measure used to move $\kap$.  Then,
in $M$ $cf\ \kap = \mu$ and $S_i = \{\bet <
\kap\mid cf \bet = \aleph_0\ {\rm and}\
o^\calF(\bet) > \bet^* + i\}$ is a stationary
subset of $\kap$ for every $i < \mu$.
Hence $S_i$ is stationary also in $V$.

\proclaim Claim 2. For every $i < \mu$
and $X \in \calP (\kap) \cap \calK (\calF)$,
$X \in \calF (\kap, \alp^* + i)$ 
iff $S_i \backslash \{ \bet < \kap |
o^\calF(\bet) > \bet^* + i$ and 
$X \cap \bet \in \calF (\bet, \bet + i) \}$
is nonstationary.

\pr Fix $i < \mu$.  $\calF(\kap, \alp^* +
i)$ is an ultrafilter over $\calP (\kap)
\cap \calK(\calF)$ so it is enough to show
that for every $X \in \calF(\kap, \alp^* +
i)$ the set $S_i \backslash \{ \bet < \kap
| o^\calF (\bet) < \bet^* + i$ and $X \cap
\bet \in \calF (\bet, \bet + i) \}$ is
nonstationary.

Suppose otherwise.  Let $X \in \calF (\kap,
\alp^* + i)$ be so that $S' = S_i
\backslash \{ \bet < \kap | o^\calF (\bet)
> \bet^* + i$ and $X \cap \bet \in \calF
(\bet, \bet + i)\}$ is stationary.
 
Without loss of generality we may assume
that $S'$ already decides the relevant
measure, i.e. for some $\gam <
o^\calF(\kap)$ $S'$ forces  the measure
$\calF(\kap, \gam)$ to be used first to move
$\kap$ in the embedding into generic
ultrapower restricted to $\calK(\calF)$.
Now, $S' \subseteq \{ \bet < \kap | o^\calF
(\bet) > \bet^* + i\}$.  So, $\gam > \gam^*
+ i$, where $\gam^*$ is the largest ordinal
$\gam^*$ below $\gam$ with $A \in
\calF(\kap, \gam^*)$.  If $\gam^* =
\alp^*$, then $\alp^* + i < \gam$ and hence
$X^* = \{ \bet < \kap | o^\calF (\bet) >
\bet^* + i$ and $X \cap \bet \in
\calF(\bet, \bet + i)\} \in \calF(\kap,
\gam)$ since this is true in the ultrapower
 of $\calK(\calF)$ by $\calF(\kap,
\gam)$.  This leads to a contradiction,
since, if $j: V \to M$ is a generic
embedding forced by $S'$, then $\kap \in
j(S')$ and $\kap \in j(X^*)$, but $S' \cap
X^* = \emptyset$. Contradiction.

If $\gam^* < \alp^*$, then also $\gam <
\alp^*$ which is impossible since there
are 
no relevant ordinals below $\alp^*$.  Also,
$\gam^*$ cannot be above $\alp^*$ since
$\alp^*$ is the last ordinal $\xi$ with $ A
\in \calF(\kap, \xi)$.\hb
$\bigsquare$

For $i < \mu$ and a set $X \subseteq \kap$
let us denote by $X^*_i$ the set $\{ \bet <
\kap | o^\calF (\bet) > \bet^* + i$ and $X
\cap \bet \in \calF(\bet, \bet + i)\}$.  By
$\Cub_\kap$ we denote the closed unbounded
filter over $\kap$ and let $\Cub_\kap \rhookup
S_i$ be its restriction to $S_i$, i.e. $\{
E \subseteq \kap | E \supseteq C \cap S_i \
\hbox{for some}\ C \in \Cub_\kap\}$.

\proclaim Claim 3.  For every $i < \mu$,
$\calF(\kap, \alp^* + i) = \{ X \in
(\calP(\kap) \cap \calK (\calF))^M | X^*_i
\in (\Cub_\kap \rhookup S_i )^M \}$.

\pr Let $X \in \calF(\kap, \alp^* + i)$,
then, by Claim 2, $X_i^* \in \Cub_\kap
\rhookup S_i$ in $V$.  But then, also in
$M$, $X^*_i \in (\Cub_\kap \rhookup S_i)^M$,
since $(\Cub_\kap)^M \supseteq
(\Cub_\kap)^V$. Now, if $X \not\in
\calF(\kap, \alp^* + i)$, then $Y = \kap
\backslash X \in \calF(\kap, \alp^* + i)$,
assuming $X \in \calP(\kap) \cap
\calK(\calF)$.  By the above, $Y^*_i \in
(\Cub_\kap \rhookup S_i)^M$.  But $X \cap Y 
= \emptyset$ implies $X^*_i \cap Y^*_i =
\emptyset$.  So $X^*_i \not\in (\Cub_\kap
\rhookup S_i)^M$.
\hb
$\bigsquare$

\proclaim Claim 4. $o^\calF(\kap) >
\alp^* + \mu$.

\pr By Claim 3, $\calF(\kap, \alp^* + i)
\in M$ for every $i < \mu$.  Hence
$(o(\kap))^M \ge \alp^* + \mu$.  But now,
in $V$, $o^\calF(\kap) \ge \alp^* + \mu +
1$.
\hb
$\bigsquare$

We actually showed more:

\proclaim Claim 5. $S \llvdash^{''}
o({\buildrul \sim \under\xi}) \ge \alp^* + \mu$ and for
every $i < \mu$ $\calF (\kap, \alp^* + i) =
\{ X \in (\calP (\kap) \cap
\calK(\calF))^{\buildrul \sim \under M} | X_i^* \in
(\Cub_\kap \rhookup S_i)^{\buildrul \sim
\under M}
\}^{''}$, where ${\buildrul \sim\under \xi}$ is a name
of the index of the first measure $\calF(\kap,
\xi)$ used to move $\kap$ and ${\buildrul
\sim \under M}$ is a generic ultrapower.

In order to complete the proof, we need to
show that every $Y \in \calF(\kap, \alp^* +
\mu)$ belongs to $\calF(\kap, \gam)$ for
unboundedly many $\gam$'s below $\alp^*$.
The conclusion of the theorem will then
follow by [Gi, Sec. 1].  So let $Y \in
\calF(\kap, \alp^* + \mu)$.  Consider a set
$Y^* = \{ \bet < \kap\mid \bet^*$  
 exists, $o^\calF(\bet) > \bet^* + \mu$  and $Y \cap \bet \in
\calF(\kap, \bet^* + \mu)\}\cup Y$.  Then $Y^* \in \cap
\{ \calF(\kap, \alp) | \alp^* + \mu \le
\alp < o^\calF (\kap)\}$.  It is enough to
show that $Y^*$ belongs to $\calF(\kap,
\gam)$ for unboundedly many $\gam$'s below
$\alp$.

\proclaim Claim 6. $S \backslash Y^*$
is nonstationary.

\pr Suppose otherwise.  Let $S' \subseteq
S\backslash  Y^*$ be a stationary set
forcing  $\calF(\kap, \xi)$ to be the
first measure used to move $\kap$ in the
ultrapower, where $\xi < o^\calF(\kap)$.
Then, by Claim 5, $\xi \ge \alp^* + \mu$.
Hence, $Y^* \in \calF(\kap, \xi)$,
which is impossible, since $Y^* \cap S'
= \emptyset$.  Contradiction.
\hb
$\bigsquare$

\proclaim Claim 7.  $\alp^*$ is a $\mu +
1$-repeat point.

\pr Let $Y^*$ be as above.  It is
enough to find $\gam < \alp^*$ such that
$Y^* \in \calF(\kap, \alp)$.  Let $C
\subseteq \kap$ be a club avoiding $S
\backslash Y^*$.  Let $N$, $\{ \tau_n |
n < \ome\}$ be as in Claim 1 (i.e. as in the
proof of 2.5.1 [Gi]) only with the club of
Claim 1 replaced by $C$ and with $Y^*
\in N$.  Then $\tau_n$'s are in $S \cap C$
and, hence in $Y^*$,  which means that for
all but finitely many $n$'s $Y^* \in
\calF(\kap, \bet^N (\tau_n))$, by [Mit
1,2], since $\tau_n$'s are indiscernibles
for $\bet^N (\tau_n)$'s.\hb
$\bigsquare$

The claim does not rule out the possibility
that some $Y^*$ reflects only boundedly
many times below $\alp^*$.  Thus, 
there is probably some $\eta< \alp^*$ such that
the  $\bet^N(\tau_n)$'s of Claim 7 are
always below $\eta$.  This means that $\bet^*_n >
\bet^N(\tau_n)$, where $\bet^*_n$ is the
stabilized value of $(\bet(\nu))^*$ for
$\nu \in d^N_n$. 
We will use Claim 5 in order to show that
this is impossible.  Namely, the following
holds:

\subheading{Claim 8}  In the notation of
Claim 7, for all but finitely many $n$'s
$\big(\bet^N (\tau_n)\big)^* = \bet^*_n$.

\pr By Claim 5, for all but nonstationary
many $\nu$'s in $S$ the following property
$(*)$ holds:\hb
$o^\calF(\nu) \ge \nu^* + \mu$ 
and for every  $i < \mu$
\enskip $\calF(\nu, \nu^* + i) = \big\{ X \in \calP
(\nu) \cap \calK(\calF) | X^*_i \in \Cub_\nu
\rhookup \{ \rho < \nu | cf \rho = \aleph_0
\ {\rm and}\ o^\calF(\rho) > \rho^* +
i\}\big\}$. 

Without loss of generality let us
assume that $(*)$ holds for every element
of $S$, otherwise just remove the
nonstationary many points. Then, preserving
notations of Claim 7, $\tau_n$'s satisfy
$(*)$.  We now show  that
ultrafilters $\calF(\tau_n, \tau_n^* + i)$
 correspond to $\calF(\kap, \bet^*_n
+ i)$ for all but finitely many $n < \ome$
and all $i < \mu$.

Let $\overline\bet^*_n$ denote
$\big(\bet^N (\tau_n)\big)^*$ and we will
drop the upper index $N$ further.  Then
$\tau_n^* + i =$ \break
$\CC$ $(\kap,
\overline\bet^*_n + i, \bet(\tau_n)
)(\tau_n)$ for every $n < \ome$, where $\CC$
is the coherence function (see [Mi1] or
[Gi]).  Suppose that $\bet^*_n \ne
\overline\bet_n^*$ for infinitely many
$n$'s.  For simplicity let us assume that
this holds for every $n < \ome$.  In the
general case only the notation is more
complicated.  There will be $X_n \in \big(
\calF(\kap, \overline\bet^*_n) \backslash
\calF(\kap, \bet^*_n) \big) \cap N$ for
every $n < \ome$, since $N$ is an
elementary submodel.  Let $n < \ome$ be
fixed.   Pick $\calK(\calF)$ - least $X_n
\in \calF(\kap, \overline\bet^*_n)
\backslash \calF(\kap, \bet^*_n)$.  Still
it is in $N$ by elementarity.  Also its
support (in the sense of [Mi1,2]) will be
below $\tau_n$, i.e. $X_n = h^N (\del)$,
for $\del < \tau_n$, where $h^N$ is the
Skolem function of $N \cap \calK(\calF)$.
The reason for this is  that $X_n$ appears
once both $\overline\bet^*_n$ and
$\bet^*_n$ appear.  But $\overline\bet_n^*$
appear below $\tau_n$ since the support of
$\tau_n$ is below $\tau_n$ and $\bet_n^*$
appear before $\tau_n$ since for $\nu \in
d_n \subseteq \tau_n$ $\big( \bet^N (\nu)
\big)^* = \bet^*_n$.  Hence $X_n \cap
\tau_n \in \calF(\tau_n, \tau^*_n)$.  Then
by
$(*)$,  $(X_n)^*_0 \in \Cub_{\tau_n}
\rhookup \{ \rho < \tau_n | cf\rho =
\aleph_0$ and $o^\calF (\rho) > \rho^*\}$.
This is clearly true also in $N$.  But
then $(X_n)_0^* \cap \cup (N \cap \tau_n)$
contains an $\ome$-club intersected with
the 
set $\{\rho < \tau_n | cf \rho = \aleph_0$
and $o^\calF(\rho) > \rho^*\}$.  Hence
$(X_n)^*_0 \cap d_n$ is unbounded in $\cup (N
\cap \tau_n)$.  Then $(X_n)^*_0 \in \calF(\kap,
\bet^*_n + i)$ for some $i$, $0 < i < \mu$,
which implies that $X_n \in \calF(\kap,
\bet_n^*)$.  Contradiction.\hb
$\bigsquare$

Combining Claims 7 and 8 we obtain that
$Y^* \in \calF(\kap, \bet^*_n + \chi)$
for some $\chi \ge \mu$, for all but
finitely many $n$'s.  Now, $\bet_n^*$'s are
unbounded in $\alp^*$ by [Gi] and hence we
have an unbounded reflection of $Y$
below $\alp^*$.
\hb
$\bigsquare$

\sect{2.~~On the strength of
saturatedness of \fakebold{$NS_\kap$}}

It was shown in [Gi] that saturatedness
of $NS_\kap$ for an inaccessible $\kap$
implies an inner model with $\exists \alp
o(\alp) = \alp^{++}$.  It was asked if 
the saturatedness of $NS_\kap^{\aleph_0}$,
i.e. the nonstationary ideal restricted to
cofinality $\ome$ already implies this.  In
this section we are going to provide an
affirmative answer.

Let us start with a ``ZFC variant" of Lemma
2.18 of [Gi].

\proclaim Proposition 2.1.  Let $V_1
\subseteq V_2$ be two models of $ZFC$.  Let
$\kap$ be a regular cardinal of $V_1$ which
changes its cofinality to $\Tet$ in $V_2$.
Suppose that in $V_1$ there is an almost
decreasing (mod nonstationary or
equivalently mod bounded) sequence of
clubs of $\kap$ of length  $(\kap^+)^{V_1}$
 so that every club of $\kap$ of
$V_1$ is almost contained in one of the
clubs of the sequence.  Assume that $V_2$
satisfies the following:
\item{(1)} $cf\ (\kap^+)^{V_1} \ge (2^\Tet)^+$ or $cf\ (\kap^+)^{V_1}
= \Tet$;
\item{(2)} $\kap > \Tet^+$.

\noindent
{\sl Then in $V_2$ there exists a cofinal
in $\kap$ sequence $\langle \tau_i \mid i <
\Tet\rangle$ consisting of ordinals of
cofinality $\ge \Tet^+$ so that every club
of $\kap$ of $V_1$ contains a final segment
of $\langle \tau_i \mid i < \Tet\rangle$.}

\subheading{Remark}  (1) If in $V_1$,
$2^\kap = \kap^+$, then clearly there exists
an almost decreasing sequence of clubs of
$\kap$ of length $\kap^+$ so that every
club of $\kap$ of $V_1$ is almost contained
in one of the clubs of the sequence.
 \hb
(2) M. Dzamonja and S. Shelah [D-Sh] using club
guessing techniques were able to replace the
condition (1) by weaker conditions.

\pr If $cf (\kap^+)^{V_1}= \Tet$ then we can simply
diagonalize over all the clubs.  So let us
concentrate on the case $cf\ (\kap^+)^{V_1}\ge
(2^\Tet)^+$.  Suppose otherwise.  Assume
for simplicity that $\Tet = \aleph_0$.  Let
$C$ be a club in $\kap$ in $V_1$.  Define
in $V_2$ a wellfounded tree $\langle T(C),
\le_C \rangle$.  Let the first level of
$T(C)$ consist of the least cofinal in
$\kap$ sequence of order type $\ome$ in
some fixed for the rest of the proof
well ordering of a larger enough portion of
$V_2$.
Suppose that $T(C) \rhookup n + 1$ is
defined.  We define Lev$_{n+1} (T(C))$.
Let $\eta \in$ Lev$_n(T(C))$.  
Let $\eta^*$ be the largest
ordinal in $T(C) \rhookup n + 1$ below
$\eta$.  We assume by induction that it
exists.
If $cf \eta
= \aleph_0$, then pick $\langle \eta_n \mid
n< \ome\rangle$ the least
cofinal sequence in $\eta$ of order type
$\ome$.  
  Let the set of immediate
successors of $\eta$, Suc$_{T(C)} (\eta)$ be
$\{\eta_n \mid n < \ome,\ \eta_n >
\eta^*\}$.

If $cf \eta \ge \aleph_1$, then consider
$\eta' = \cup (C \cap \eta)$.  If $\eta' =
\eta$, then let Suc$_{T(C)} (\eta) =
\emptyset$.  If $\eta^*< \eta' < \eta$, then
let Suc$_{T(C)} (\eta) = \{ \eta'\}$.
Finally, if $\eta' \le \eta^*$ then let
Suc$_{T(C)} (\eta) = \emptyset$.  This
completes the inductive definition of
$\langle T(C)$, $\le_C\rangle $.  Obviously, it is 
wellfounded and countable.  Let $T^* (C)$
denote the set of all endpoints of $T(C)$
which are in $C$.  Notice, that by the
construction any such point is of
uncountable cofinality.  Also, $T^*(C)$ is
unbounded in $\kap$, since otp$(C) = \kap$
and $\kap > \aleph_1$.

There must be a club $C_1 \subseteq C$ in
$V_1$  avoiding unboundedly many points of
$T^* (C)$, since otherwise the sequence
$\langle\tau_i | i < \aleph_0\rangle$ required by the
proposition could be taken from $T^*(C)$.
 This means,in
particular, that for every $\alp < \kap$
there will be 
$$\overline \nu = \langle
\nu_1 \nek \nu_n\rangle \in T(C) \cap
T(C_1)$$
so that
\item{(a)} $cf\ \nu_n > \aleph_0$;
\item{(b)} Suc$_{T(C)} (\nu_n) = \{\nu_{n+1}\}$
for some $\nu_{n+1} \in C \backslash \alp$;
\item{(c)} either
\itemitem{(c1)} Suc$_{T(C_1)}
(\nu_n) = \emptyset$
 
\noindent
or
 \itemitem{(c2)} for some
$\rho \in (C_1 \cap \nu_{n+1} ) \backslash
\alp$\quad  Suc$_{T(C_1)} (\nu_n) = \{ \rho\}$.

Now define a sequence $\langle C_\alp \mid
\alp < (2^{\aleph_0} )^+ \rangle$ of clubs
so that 
\item{(1)} $C_\alp$ is a club in $\kap$ in
$V_1$;
\item{(2)} if $\bet < \alp$ then $C_\alp
\backslash C_\bet$ is bounded in
$\kap$;
\item{(3)} $C_{\alp + 1}$ avoids
unboundedly many points of $T^* (C_\alp)$.

Since $cf\
(\kap^+)^{V_1} \ge (2^{\aleph_0})^+$ and in
$V_1$ there is an almost decreasing (mod
bounded) sequence of $\kap^+$-clubs
 generating
the club filter, there is no problem in
carrying out the construction of $\langle
C_\alp \mid \alp < (2^{\aleph_0})^+
\rangle$ satisfying (1)--(3).

Shrinking the set of $\alp$'s if necessary
we can assume that for every $\alp, \bet <
(2^{\aleph_0})^+$ $\langle T(C_\alp),\
\le_{C_2}, \le \rangle$ and $\langle
T(C_\bet),\ \le_{C_\bet}, \le \rangle$ are
isomorphic as  trees with ordered levels.

Let $\langle \kap_m \mid m < \ome\rangle$
be the least cofinal in  $\kap$ sequence.

Let $\alp < \bet < (2^{\aleph_0})^+$.
Since $C_\bet$ is almost contained in
$C_{\alp + 1}$, it avoids unboundedly many
points in $T^*(C_\alp)$.  So for every $m <
\ome$ there is $\overline\nu = \langle
\nu_1 \nek \nu_n \rangle \in T(C_\alp) \cap
T(C_\bet)$ so that
\item{(a)} $cf\ \nu_n > \aleph_0$;
\item{(b)} Suc$_{T(C_\alp)} (\nu_n) =
\{\nu^\alp_{n+1}\}$ for some $\nu_{n+1}^\alp
\in C_\alp \backslash \kap_m$;
\item{(c)} 
 for some $\nu_{n+1}^\bet
\in (C_\bet \cap \nu_{n+1}^\alp) \backslash
\kap_m\quad {\rm Suc}_{T(C_\bet)}(\nu_n) =
\{\nu_{n+1}^\bet\}$.

Thus, pick $\ell > m$ so that $C_\bet
\backslash \kap_{\ell - 1} \subseteq
C_\alp$.  We consider subtrees
$$T(C_\gam)_\ell = \{ \overline \eta \in
T(C_\gam) | \exists k \ge \ell \ \overline
\eta\ \ge_{C_\gam}  \langle \kap_k\rangle\}$$
where $\gam = \alp, \bet$.

Let $\pi$ be an isomorphism between
$T(C_\alp)$ and $T(C_\bet)$ respecting
the order of the levels.  Notice, that the
first level in both trees is the same
$\{\kap_i | i < \ome\}$.  Hence, $\pi$ will
move $T(C_\alp)_\ell$ onto
$T(C_\bet)_\ell$.

Pick the maximal $n < \ome$ such that $\pi
$ is an identity on $\big(T(C_\alp)_\ell
\big)\rhookup n+1$.  It exists since
$T^*(C_\alp) \backslash C_\bet$ is
unbounded in $\kap$.  Now let $\nu$ be the
least ordinal in ${\rm Lev}_{n + 1}
\big(T(C_\alp)_\ell\big)$ such that $\pi (\langle
\nu_1 \nek \nu_n, \nu\rangle ) \ne \langle \nu_1 \nek
\nu_n, \nu\rangle$, where $\langle \nu_1 \nek \nu_n\rangle$
is the branch of $T(C_\alp)_\ell$ leading
to $\nu$.  

Consider $\nu_n$.  If $cf \nu_n =
\aleph_0$, then we are supposed to pick the
least cofinal in $\nu_n$ sequence $\langle
\nu_{ni} | i < \ome \rangle$ and the maximal
element $\nu_n^*$ of the tree $T(C_\alp)$
below $\nu_n$.
Suc$_{T(C_\alp)}  (\nu_n)$ will be $\{
\nu_{ni} | i < \ome$ and $\nu_{ni} >
\nu_n^*\}$.  Notice that $\nu_n^* \ge
\kap_{n-1}$ by the definition of the tree
$T(C_\alp)$.  Hence, either $\nu^*_n =
\kap_{n-1}$ or $\nu_n^* \in T(C_\alp)_\ell
\rhookup n + 1$ since elements of $T(C_\alp)$
which are above $\kap_{n-1}$ in the tree
order are below it as ordinals.  But since
$T(C_\alp)_\ell \rhookup n + 1 = T
(C_\bet)_\ell \rhookup n + 1$ and
$\kap_{\ell - 1} \in T(C_\bet)$, the same
is true about Suc$_{T(C_\bet)} (\nu_n)$,
i.e. it is $\{ \nu_{ni} | i < \ome $ and
$\nu_{ni} > \nu^*_n\}$.  Then $\pi$ will be
an identity on Suc$_{T(C_\alp)} (\nu_n)$ and,
in particular, will not move $\nu$.
Contradiction.

So $cf \nu_n$ should be above $\aleph_0$.
Once again the maximal elements of
$T(C_\alp) \rhookup n + 1$ and $T(C_\bet)
\rhookup n + 1$ below $\nu_n$ are the same.
Let $\nu^*_n$ denote this element.  Now,
$\nu \in {\rm Suc}_{T(C_\alp)} (\nu_n)$,
hence $\nu = \cup (C_\alp \cap \nu_n)\enskip
 \nu_n^* <
\nu < \nu_n$ and Suc$_{T(C_\alp)} (\nu_n) =
\{ \nu\}$ by the definition of the tree
$T(C_\alp)$.  $\pi$ is an isomorphism, so
Suc$_{T(C_\bet)} (\nu_n) \ne \emptyset$.
By the definition of the tree $T(C_\bet$),
 $\nu^* < \nu' < \nu_n$ and
Suc$_{T(C_\bet)} (\nu_n) = \{ \nu'\}$ where
$\nu' = \cup (C_\bet \cap \nu_n)$.  By the
choice of $\nu$, $\nu \ne \nu'$.  But $\nu,
\nu' > \kap_{\ell - 1}$ and $C_\bet
\backslash \kap_{\ell - 1} \subseteq
C_\alp$, so $\nu' \in C_\alp$.  Hence $\nu'
< \nu$ and the sequence $\langle\nu_1 \nek
\nu_n\rangle$ is as desired.

  Let $\langle T,
\le_T, \le \rangle$ be a countable tree
consisting of countable ordinals with the
usual order $\le$ between them isomorphic
to $\langle T(C_\alp), \le_{C_\alp}, \le
\rangle$ $\bigl(\alp < (2^{\aleph_0})^+\bigr)$.
Define a function $h:[(2^{\aleph_0})^+]^2
\to \ome$ as follows:\hb
$f(\alp, \bet) =$ the minimal element of
$T$ corresponding to some
$$\overline \nu = \langle \nu_1 \nek
\nu_n\rangle \in T(C_\alp) \cap T(C_\bet)$$
satisfying the conditions (a), (b) and
(c).

By Erd\" os--Rado there exists a homogeneous
infinite set $A \subseteq (2^{\aleph_0})^+$.
Let $\langle \alp_n \mid n < \ome\rangle$
be an increasing sequence from $A$.  Then
there is $\overline\nu = \langle \nu_1\nek
\nu_n\rangle \in \bigcap\limits_{m < \ome}
T(C_{\alp_m})$ witnessing (a), (b), (c).
But by (c), $\nu_{n+1}^{\alp_m} >
\nu_{n+1}^{\alp_m+1}$ for every $m < \ome$.
Contradiction.\hb
 $\bigsquare$

Suppose now that there is no inner model of
$\exists \alp o(\alp) = \alp^{++}$.
The following follows easily from
Proposition 2.1 and the Mitchell Covering
Lemma [Mi3].

\proclaim Proposition 2.2.  The final
segment of the sequence
$\langle \tau_n\mid n < \ome\rangle$
consists of indiscernibles for
$\kap$.

\pr Suppose otherwise.  Then by the
Mitchell Covering Lemma [Mi3] there is $h
\in \calK(\calF)$ and $\del_n < \tau_n$ $(n
< \ome)$ such that $h(\del_n) \ge \tau_n$
for infinitely many $n$'s.  Define a club
in $\calK (\calF)$:
$$C = \{ \nu < \kap | h^{''} (\nu) \subseteq
\nu\}\ .$$
Then, by the choice of $\langle\tau_n | n <
\ome\rangle$, there is $n_0 < \ome$ such that for
every $n \ge n_0$ $\tau_n \in C$,  which is
impossible.  Contraction.
\hb
$\bigsquare$

Let us show that $o(\kap) =
\kap^{++}$ if $NS^\Tet_\kap$ is saturated.

\proclaim Proposition 2.3. \footnote{$^*$}{\smalltype 
This result was recently improved by S.
Shelah and the author to $0 = 1$.}
 Suppose that
$NS_\kap^\Tet$ is saturated over an
inaccessible $\kap$, then $o(\kap) =
\kap^{++}$.

\subheading{Remark} 
By Mitchell [Mi1] it follows for
successor $\kap$'s and moreover, by Shelah
[Sh] it is impossible for successor
cardinal $\kap$ which is  above $\Tet^+$.

\pr Let for simplicity $\Tet = \aleph_0$. Suppose that
$o(\kap) < \kap^{++}$.  We call an ordinal
$\alp$ a relevant ordinal if some $S \in
(NS_\kap^{\aleph_0})^+$ forces the
measure $\calF(\kap, \alp)$ (of the core
model $\calK(\calF)$) to be used as the first measure
moving $\kap$ in the generic embedding
restricted to $\calK(\calF)$.  By Mitchell
[Mi1,2], such a restriction is an iterated
ultrapower of $\calK(\calF)$.  Let us call the
corresponding measure $\calF(\kap, \alp)$
-- a relevant measure.

Since $NS_\kap^{\aleph_0}$ is saturated,
the total number of relevant measures is at
most $\kap$.  Let $\langle A_\alp \mid \alp
< \chi \le \kap\rangle$ be a maximal antichain
such that $A_\alp$ forces ``$\alp$ to be a
relevant ordinal".  Without loss of
generality $A_\alp$'s are pairwise disjoint
and min$A_\alp  > \alp$.  Also it is
possible to pick each $A_\alp$ in
$\calF(\kap, \alp)$ using $o(\kap) <
\kap^{++}$, but it is not needed for the
rest.  Let us assume that $\chi = \kap$.  The
case $\chi < \kap$ is similar and even slightly
simpler.

The set $A = \bigcup\limits_{\alp < \kap}
A_\alp$ contains an $\aleph_0$--club.
Since $NS_\kap^{\aleph_0}$ is precipitous,
we can assume that every $\nu \in A$ is
regular in $\calK(\calF)$.  Otherwise just remove
from $A$ nonstationary many $\nu$'s of
cofinality $\aleph_0$ which are singular in
$\calK(\calF)$.

Let $\alp_{\min}$ be the least relevant
ordinal.  Form a generic ultrapower $M$
with $A_{\alp_{\min}}$ in a generic
ultrafilter $G \subseteq
(NS_\kap^{\aleph_0})^+$.  Applying
Proposition 2.1 to $V$ and $V[G]$, we find
a cofinal in $\kap$ sequence $\langle
\tau_n\mid n < \ome\rangle$ consisting of
ordinals of uncountable cofinality such
that every club of $\kap$ of $V$ contains
its final segment.  Since
$NS_\kap^{\aleph_0}$ is saturated in $V$,
$\langle \tau_n \mid n < \ome\rangle \in
M$.  We are going to use $\langle \tau_n
\mid n < \ome\rangle$ in order to recover
$\calF(\kap, \alp_{\min})$ inside $M$,  which
is impossible, since it is used already to
move $\kap$ and hence cannot be in
$\calK(\calF)^M$.

We proceed as follows.  For every $\alp \in
A \backslash \{\alp_{\min}\}$ pick in $V$ a
function $g_\alp \in {}^\kap\kap$ which is
forced by $A_\alp$ to represent $\alp_{\min}$
in a generic ultrapower.  By saturatedness
it is possible.  Under $o(\kap) <
\kap^{++}$ we can find such $g_\alp$ in
$\calK(\calF)$.

For a set $X \in \calF(\kap, \alp_{\min})$ we
define in $V$ a set\hb
$C_X = \Big\{\nu<\kap \mid cf\ \nu \ne
\aleph_0\ {\rm or}\ \bigl(cf\ \nu =
\aleph_0$
and then either $\nu \in A_{\alp_{\min}} \cap
X$ or $\nu \in A_\alp$ for some $\alp \in A
\backslash \{\alp_{\min}\}$ and then $X \cap
\nu \in \calF(\nu, g_\alp (\nu)\bigr)\Big\}$.

Then $C_X$ contains a club.   $C_X$
is in $M$ as well as $\langle A_\alp \mid \alp
\in A\rangle$ and $\langle g_\alp \mid \alp
\in A\rangle$.  Moreover it is has 
the same definition as in $V$.

But now, in $M$, we may define a set $D = \{
X \subseteq \kap \mid X \in (\calK)^M$,
$C_X$ contains a final segment of $\langle
\tau_n \mid n < \ome\rangle \}$.  Then $D
\supseteq \calF (\kap, \alp_{\min})$, since
for every $X \in \calF(\kap, \alp_{\min})$ $X
\in M$ and $C_X$ is a club of $V$.  On the
other hand, if $X \in D$, then $X \in
\calK(\calF)^M \cap \calP(\kap) \subseteq \calK(\calF)
\cap \calP(\kap)$.  If $X \not\in
\calF (\kap, \alp_{\min})$, then $Y = \kap
\backslash X$ belongs to $\calF{(\kap,
\alp_{\min})}$ and hence $C_Y$ contains a
final segment of $\langle \tau_n \mid n <
\ome\rangle$,  as does  $C_X$.  But
this is impossible.  Thus let $C'_X, C'_Y$
be the clubs of limit points of $C_X$ and
$C_Y$.  Let $\tau_n \in C'_X \cap C'_Y$,
then there is some $\nu < \tau_n$, $cf \
\nu = \aleph_0$ and $\nu \in C_X \cap C_Y$
since $cf{\tau_n} > \aleph_0$.
But then for some unique $\alp \in A$, $\nu
\in A_\alp$ which implies $\nu \in X \cap
Y$, in the case $\nu \in A_{\alp_{\min}}$,
 or $X
\cap \nu$, $Y \cap \nu$ in $\calF(\nu,
g_\alp(\nu))$ otherwise.  Which is
impossible since $Y$ and $X$ are disjoint.

So $D = \calF(\kap, \alp_{\min})$.  Hence
$\calF(\kap, \alp_{\min})\in M$.
Contradiction.\hb
$\bigsquare$

We think that the methods of [Gi-Mi] can be used in order
to push the strength of ``$NS_\kap$ saturated" to a strong cardinal.

\sect{3.~~On the strength of
precipitousness of a nonstationary ideal}
\vskip -.3truecm
\ce{\bf over
an inaccessible}
\medskip

We are going to show that the assumptions
used in [Gi] making $NS_\kap$ precipitous
$\bigl( ( \ome,\kap^+ + 1)$ - repeat point
$\bigr)$ and $NS_\kap^{\aleph_0}$
precipitous $\bigl( (\ome, \kap^+)$--repeat
point$\bigr)$ over an inaccessible $\kap$
can be weakened to an $(\ome, \kap + 1)$--repeat
point and to an $(\ome, \kap)$-repeat point,
respectively.  This is quite close to the 
equiconsistency, since by [Gi], an $(\ome, <
\kap)$--repeat point is needed for the
existence of such ideals. 

\proclaim Theorem 3.1.  Suppose that there
exists an $(\ome, \kap + 1)$--repeat point
over $\kap$.  Then in a generic extension
preserving inaccessibility of $\kap$,
$NS_\kap$ is a precipitous ideal.

The proof combines constructions of [Gi]
and [Gi1].  We will stress only the new
points.

\medskip

\noindent
{\bf Sketch of the Proof:}
Let $\alp < o(\kap)$ be an $(\ome,
\kap+1)$--repeat point for $\langle
\calF(\kap, \alp')\mid \alp' <
o(\kap)\rangle$,  i.e. $cf\ \alp =
\aleph_0$ and for every $A \in \cap \{
\calF(\kap, \alp^* + i)\mid i \le \kap\}$
there are unboundedly many $\bet$'s in $\alp$ 
 such that $\bet + \kap < \alp$
and $A \in \cap \{ \calF(\kap, \bet+i) \mid
i \le \kap\}$.

As in [Gi] we first define the iteration
$\calP_\del$ for $\del$ in the closure of
$\{ \bet \le \kap \mid \bet$ is an
inaccessible or $\bet = \gam + 1$ for an
inaccessible $\gam\}$.  On limit stages as
in [Gi] the limit of [Gi2] is used.  Define
$\calP_{\del + 1}$.  If $o(\del) \ne \bet +
\del$ or $o(\del) \ne \bet + \del + 1$ for
some $\bet$ then $\calP_{\del + 1} =
\calP_\del * C(\del^+) * \calP(\del,
o(\del))$ exactly as in [Gi], where
$C(\del^+)$ is the Cohen forcing for adding
$\del^+$ functions from $\del$ to $\del$
and $\calP(\del, o(\del))$ is a forcing
used in [Gi] for changing cofinalities
without adding new bounded sets.

Now let $o(\del) = \bet + \del$ for some
ordinal $\bet, \bet > \del$.  First we force as above
with $C(\del^+)$.  
\subheading{Case 1} The value of the
first Cohen function added by $C(\del^+)$
on $0$ is not $0$. \hb
Then we force as above with $\calP(\del, o(\del))$.  
\subheading{Case 2} 
The value of the first Cohen function added
by $C(\del^+)$ on 0 is 0.  \hb
 Then we are going to shoot a
club through $\cap\{\calF(\del, \bet + i)
\mid i < \del\}$ using the forcing notion
$Q$ described below.

$Q = \{ \langle c, e\rangle \mid c
\subseteq \del\ {\rm closed}, |c| < \del,
\ e \subseteq \cap\{\calF (\del, \bet + i)
\mid i < \del\}, \ |e| < \del\}$\hb
 $\langle
c_1, e_1\rangle \le \langle c_2,
e_2\rangle$ iff $c_2$ is an end-extension of
$c_1$, $e_1 \subseteq e_2$ and for every $A
\in e_1$, $c_2 \backslash c_1 \subseteq A$.
Now every regular $i < \del$ forcing
with $\calP(\del, \bet + i)$ produces a
club through $\cap\{\calF(\del, \bet + j)
\mid j < i\}$ changing cofinality of $\del$
to $i$.  
Thus $Q$ contains an $i$-closed dense
subset in any $\calP(\alp, \bet +
i)$-generic extension of $V^{\calP_\alp * C
(\alp^+)}$. Based on
this observation, we are going to use here
the method of [Gi1]. It makes the iteration
of such forcings $Q$ possible.

If $o(\del) = \bet + \del + 1$ for some
$\bet,\ \bet > \del$, then we  combine both
previous cases together inside the Prikry
sequence produced at this stage.

Namely, we proceed as follows.  Let $i:V
\to M \simeq Ult \bigl(V, \calF(\del, \bet +
\del)\bigr)$.  We consider also the second
ultrapower, i.e. $N \simeq Ult \bigl(M, \calF(i(\del),
i(\bet) + i(\del)\bigr)$.  Let $k:M \to N$
and $j = k \circ i:V \to N$ be the
corresponding elementary embeddings.  Then,
in $N$, $o(\del) = \bet + \del$ and
$o(i(\del)) = i(\bet) + i(\del)$.  So, in
$N$, both $\del$ and $i(\del)$ are of the
type of the previous cases.  We want to
deal with $\del$ as in Case 1 and with
$i(\del)$ as in Case 2.  This can easily be
arranged, since we are free to change one
value of a Cohen function responsible for
the switch between Cases 1 and 2.  The next
stage will be to define an extension
$\calF^* (\del, \bet + \del)$ of
$\calF(\del, \bet + \del) \times \calF
(\del, \bet + \del)$ in $V[G_\del]$, where
$G_\del \subseteq \calP_\del$ is generic.
For this use [Gi1] where $N$ was
first stretched by using the direct limit of
$\langle \calF\bigl(i(\del), i(\bet) +
\xi\bigr) \mid \xi < i(\del)\rangle$.
Finally we force a Prikry sequence using
$\calF^*(\del, \bet + \del)$.  Notice that
the following holds:\hb
{\it $(*)$ if $\langle \langle
\del_n, \rho_n \rangle \mid n < \ome
\rangle$ is such a sequence then both
$\langle \del_n \mid n < \ome \rangle$ and
$\langle \rho_n \mid n < \ome\rangle$ are
almost contained in every club of $\del$ of
$V$.}
\hb
 Simply because $\langle \del, i(\del)
\rangle \in j(C)$ for a club $C \subseteq
\del$ in $V$.

This completes the definition of
$\calP_{\del + 1}$ and hence also the
definition of the iteration.

The intuition behind this is as follows.
We  add a club subset to every set
$A \in \bigcap \{\calF (\kap, \alp + i) |
i \le \kap\}$. $\alp$ is $(\ome, \kap +
1)$-repeat point, so $A$ reflects
unboundedly many times  in $\alp$, i.e.
$A \in \cap \{ \calF(\kap, \bet + i)| i \le
\kap\}$ for unboundedly many $\bet$'s in
$\alp$.

Reflecting this below $\kap$, we will have
$A \cap \del \in \cap \{ \calF(\del, \gam +
1) | i \le \del \}$, where $o(\del) = \gam
+ \del$.  In [Gi, Sec. 3], we had $(\alp,
\kap^+ + 1)$-repeat point which translates
to $\cap \{ \calF(\del, \gam+i) | i \le
\del^+\}$.  Then just the forcing
$\calP(\del, o(\del))$ will add a club
through every set in $\cap \{ \calF(\del,
\gam + i) | i \le \del^+\}$.  Here our
assumptions are weaker and we use the
forcing $Q$ instead.  There are basically
two problems with this: iteration and
integration with $\calP(\del, \bet)$'s.  For
the first problem the method of [Gi1] is
used directly.  The problematic point with
the second is that once using $Q$ we break
the Rudin-Keisler ordering of extensions of
$\calF(\del, \bet)$'s used in $\calP(\del,
o(\del))$.  In order to overcome this
difficulty, we split the case $o(\del) =
\bet + \del$ into two.  Thus in Case 1 we
keep Rudin-Keisler ordering and in Case 2
force with $Q$.  Finally at stages $\alp$
with $o(\del) = \bet + \del + 1$ both cases
are combined in the fashion described
above.
The rest of the proof is as in [Gi, Sec. 3].

The following obvious changes needed to be
made: instead of $E \in \cap \{ \calF(\kap,
\bet) | \alp < \bet \le \alp + \kap^+\}$ we
now deal with $E \in \cap \{ \calF(\kap,
\bet)| \alp < \bet \le \alp + \kap\}$ and
instead of $E(\kap^+)$ there we use
$E(\kap) = \{ \del \in E|$ there is
$\overline\del$ s.t. $o^\calF(\del) =
\overline\del + \kap$ and $\del \cap E \in
\cap \{ \calF(\del, \del') | \overline\del
\le \del' < \overline \del + \kap^+\}$
which belongs to $\calF(\kap, \bet + \kap)$
for unboundedly many $\bet$'s in $\alp$.
Lemmas 3.2-3.5 of [Gi] have the same proof
in the present context.  The changes in the
proof of Lemma 3.6 of [Gi] (actually the
claim there) use the method of iteration of
$Q$'s and the principal $(*)$.

If we are not concerned about a regular cardinal,
then the same construction starting with
an $(\ome, \kap)$--repeat point turns
$NS_\kap^{Sing}$ into a precipitous ideal.
So the following holds:

\proclaim Theorem 3.2.  Suppose that there
exists an $(\ome, \kap)$--repeat point over
$\kap$.  Then in a generic extension
preserving inaccessibility of $\kap$,
$NS_\kap^{Sing}$ is a precipitous ideal.

\sect{4.~~Open problems}

\subheading{A.~~Saturatedness}

Are the following statements consistent:

\noindent
1$^*$.~~$NS_\kap$ saturated over
an inaccessible $\kap$.

\noindent
2\footnote{$^*$}{\smalltype No by a recent
result of S. Shelah and the author.}.~~$NS_\kap^\tet$
saturated over an inaccessible $\kap$ for
some fixed cofinality $\tet$.

\noindent
3.~~~$NS_{\kap^+}^\kap$
saturated for a cardinal $\kap \ge
\aleph_1$.

Known that $NS_{\aleph_1}$ can be saturated
[St-Van-We], [Fo-Ma-Sh], [Sh-Wo],
$NS_{\kap^+}^\tet$ cannot for $\tet < \kap$
[Sh].  By [Je-Wo] $NS_\kap \rhookup Reg$
can be saturated over inaccessible.

Test question:

\noindent
4.~~~Let $E \subseteq Reg \cap \kap$, $\kap$
inaccessible, there is no sharp for a
strong cardinal.  Suppose $NS_\kap \rhookup
E$ is precipitous (or saturated), $E
\|\!\!\!-\!\!\!- U$ is the normal measure of
the extender used to move $\kap$ by a
generic embedding restricted to the core
model.  

Is there a nonstationary set in
$U$ in a generic ultrapower?

\subheading{B.~~Precipitous ideals}

\noindent
1.~~Is the strength of
$NS_\kap^{\aleph_0}$ precipitous over
an inaccessible $\kap$ $(\ome, <\kap)$-repeat
point?

\noindent
2.~~Can a model for $NS_\kap$ precipitous
over an inaccessible $\kap$ be constructed
from something weaker than an $(\ome, \kap +
1)$-repeat point?

\noindent
3.~~What is the strength of $NS_\kap$
precipitous over the first inaccessible?

The upper bound for (3) is a Woodin
cardinal, see [Sh-Wo].  If it is possible to
construct a model with $NS_\kap^{\aleph_0}$
precipitous from an $(\ome, <\kap)$-repeat
point, then we think that this assumption
is also sufficient for (2) and (3).

A test question:
\hb
4.~~How strong is ``there is a precipitous
ideal over the first inaccessible"?

By [Sh-Wo] a Woodin cardinal suffices.  On
the other hand one can show that at least
$o(\kap) = \kap$ is needed.

\vskip 1truecm

\references{60
}

\ref{[D-Sh]} M.Dzamonja and S.Shelah,Some results on squares,outside
guessing clubs.

\ref{[Fo-Ma-Sh]} M. Foremen, M. Magidor and
S. Shelah, Martin Maximum I, {\it Ann. of
Math.} {\bf 127} (1988) 1-47.

\ref{[Gi]} M. Gitik, Some results on
nonstationary ideal, {\it Israel J. of
Math.}

\ref{[Gi1]} M. Gitik, On clubs consisting
of former regulars, to appear in {\it J. of
Sym. Logic}.

\ref{[Gi2]} M. Gitik, Changing cofinalities
and the nonstationary ideal, {\it Israel J.
of Math.} {\bf 48} (1984) 257-288.

\ref{[Gi-Mi]} M. Gitik and W. Mitchell,
Indiscernible sequences for extenders and the singular cardinal hypothesis.

\ref{[Je-Wo]} T. Jech and H. Woodin,
Saturation of the closed unbounded filter,
{\it Trans. Am. Math. Soc.} {\bf 292(1)}
(1985) 345-356.

\ref{[Mi1]} W. Mitchell, The core model for
sequences of measures I,
{\it Math. Proc.
of Cambridge Math. Soc.} {\bf 95} (1984)
41-58.

\ref{[Mi2]} W. Mitchell, The core model for
sequences of measures II. 

\ref{[Mi3]} W. Mitchell, Applications of
the Core Model for sequences of
measures, {\it Trans. Am. Math. Soc.}
{\bf 299} (1987), 41-58.

\ref{[Sh]} S. Shelah, Proper Forcing,
Springer Verlag, 1982.

\ref{[Sh-Wo]} S. Shelah and H. Woodin.

\ref{[St-Van-We]} J. Steel and R. Van
Wesep, Two consequences of determinancy
consistent with choice, {\it Trans. Am.
Math. Soc.} {\bf 272(1)} (1982) 67-85.
\end